\documentclass{article}

\usepackage{graphicx}

\usepackage{placeins}

\usepackage{tikz}
\usetikzlibrary{shapes.geometric, arrows.meta, positioning}
\usetikzlibrary{positioning}

 \tikzstyle{startstop} = [rectangle, rounded corners, minimum width=3cm, minimum height=1cm,text centered, draw=black, fill=blue!20]
\tikzstyle{process} = [rectangle, minimum width=3cm, minimum height=1cm, text centered, draw=black, fill=gray!20]
\tikzstyle{decision} = [diamond, aspect=2, text centered, draw=black, fill=orange!20]
\tikzstyle{arrow} = [thick,->,>=stealth]

\usepackage[applemac]{inputenc}

\newcommand{\tras}{^{\mbox{\scriptsize tr}}}
\newcommand{\E}{\mbox{\bf E}}
\newcommand{\var}{\mbox{\bf Var}}
\newcommand{\cov}{\mbox{\bf Cov}}
\newcommand{\N}{\mbox{\bf N}}

\newcommand{\R}{\mbox{\bf R}}

\newcommand{\pr}{\mbox{\bf P}}
\newcommand{\un}{\mbox{\bf 1}}
\newcommand{\qed}{\hfill\framebox{$ $}}
\newcommand{\ru}{\hspace{.1em}\rule[-.5ex]{.15em}{2.2ex}\hspace{.1em}}
\newcommand{\sN}{\mbox{\scriptsize N}}
\newcommand{\bm}[1]{\mbox{\boldmath $#1$}}
\newcommand{\sbm}[1]{\mbox{\scriptsize\boldmath $#1$}}

\newtheorem{theorem}{Theorem}
\newtheorem{lemma}{Lemma}
\newtheorem{corollary}{Corollary}
\newtheorem{remark}{Remark}
\newtheorem{definition}{Definition}

\begin{document}

\title{Brownian sheet and uniformity tests on the hypercube}

\author{Alejandra Caba\~na\footnote{Corresponding Author. E-mail: AnaAlejandra.Cabana@uab.cat} and Enrique M. Caba\~na}

\maketitle

\begin{abstract}
A construction of $p$-parameter Brownian sheet on the hypercube $C=[0,1]^p$ as a sum of $2^p$ independent Gaussian processes is obtained. The terms are closely related to Brownian pillows, and the probability laws of their $L^2(C)$ squared norms are computed.
This allows us to propose consistent tests of uniformity for samples of i.i.d. random vectors on $C$. 
A comparison of powers of the new tests with those of several uniformity 
tests found in the statistical literature completes the article.


\end{abstract}

\section{Introduction}\label{sec:intro}

Testing for uniformity in the hypercube $C:=[0,1]^p$
 (where  $p\geq1)$, in spite of the extreme particularity of such probabilistic model,
is interesting for several reasons, specially in areas such as statistics, machine learning, and computational mathematics.

Many Monte Carlo methods generate samples that should ideally be uniformly distributed over $[0,1]^p$.
Testing uniformity ensures that these methods produce correct and unbiased samples, which is critical for accurate integration, optimization, or probabilistic modeling.

Uniformity in the hypercube may arise as a consequence of transforming a multivariate distribution. 
For example, probability integral transforms via cumulative distribution functions (see \cite{rosenblatt1952}) map arbitrary distributions 
into the uniform distribution on $[0,1]^p$. Testing for uniformity can validate the correctness of these transformations and evaluate model fit.
Even testing composite hypothesis as goodness--of--fit to multivariate normality on $\mbox{R}^p$ can be reduced via standardization and probability integral transformation to assess uniformity on $[0,1]^p$.

Several authors base their uniformity tests on quite different statistics, namely, the distribution of distances between the data or the use of normal quantiles (\cite{Yang2015}), the depth of the elements of the sample (\cite{Hegazy1975}), their distance--to--boundary (\cite{berrendero2006testing}), random graphs over points of the sample (\cite{ebner2020testing}) or minimal covering trees of the sample graph (\cite{PETRIE2013253}). Some of these tests are included in the R (The R Project for Statistical Computing) package SHT.

.

It is striking that the use of sophisticated statistical methods is preferred over basing decisions solely on the distance between the empirical distribution of the sample and the uniform distribution. This preference is attributable to the fact that the distributions of the empirical processes under the null hypothesis of uniformity not only depend on the dimension but also on the sample size. Moreover, neither these distributions nor their asymptotic forms for large samples are well understood. 
We base our tests in a decomposition of the Brownian sheet into a sum of independent processes, that generalizes 
the well known decomposition $w(t)=tw(1)+b(t)$ of a standard Wiener process $w(t)$ on $[0,1]$ as the sum of a Brownian bridge $b(t)$ and the ``ramp" $tw(1)$.

The result is not surprising but we include a straightforward presentation in Section \ref{desco} for the sake of completeness. 
This decomposition will enable us to express the empirical process  as a sum of asymptotically independent processes. 
This key insight allows us to propose tests similar in spirit to the Cram\'er--von Mises test but with fewer drawbacks than a direct generalization would entail.

All the power comparisons in Section \ref{compar} are made via a Monte Carlo approximation of the distribution of the test statistics, since the asymptotic distribution is described by a series of random terms that is not suitable for calculation, and leads to conservative tests.

\bigskip 

\section{Decomposition of $p$-Brownian sheet as a sum of $2^p$ independent Gaussian terms}\label{desco}
The probabilistic result stated in Section \ref{Brownianramps} is based on a general property of functions with domain $C$ that vanish on the set $\partial^-C$ of points $\bm t=(t_1,t_2,\dots,t_p)\tras\in C$ with at least one component equal to zero.  These functions have a unique decomposition as a sum of terms belonging to a special class that we call {\em ramps}. 

\subsection{Ramp components of a function $g:C\rightarrow \R$ that vanishes on $\partial^-C$.}

Given a subset $H$ of the set $J=\{1,2,\dots,p\}$ and a point $\bm t=(t_1,t_2,\dots,t_p)\in C$, let us denote $\bm t_H$ the point with coordinates $(\bm t_H)_j=t_j$ if $j\in H$ and $(\bm t_H)_j=1$ if $j\not\in H$. Then introduce the sets $C_H=\{\bm t_H:\bm t\in C\}$ and $\partial_H C_H=\{\bm t_H\in C_H: \mbox{for some }j\in H, t_j=0 \mbox{ or }t_j=1\}$. Each $C_H$ is said to be a {\em face} of $C$ of dimension equal to the cardinal $\#H$ of $H$. 
The face $C_{\emptyset}$ has only the element $\bm 1$ with all coordinates equal to one, and $C_J$ is equal to $C$.

A function $g_H:C_H\rightarrow\R$ that vanishes on $\partial_H C_H$ shall be called a {\em tent}, or an {\em $H$-tent} if the domain is wanted to be mentioned explicitely.

To each $H$-tent $g_H$ we associate an extension to the domain $C$, namely $R_{g_H}(\bm t)=\ru\bm t_{J\setminus H}\ru g_H(\bm t_H)$, where $\ru(t_1,t_2,\dots,t_p)\ru$ denotes the product $\prod_{j=1}^pt_j$. Such extension is said to be an {\em $H$-ramp}. 
More precisely, the $H$-ramp associated to the $H$-tent $g_H$.
\begin{figure}
\vspace{-5mm}
\caption{Tents and ramps in $[0,1]^2$}
\vspace{-15mm}
\centerline{\includegraphics[scale=0.55]{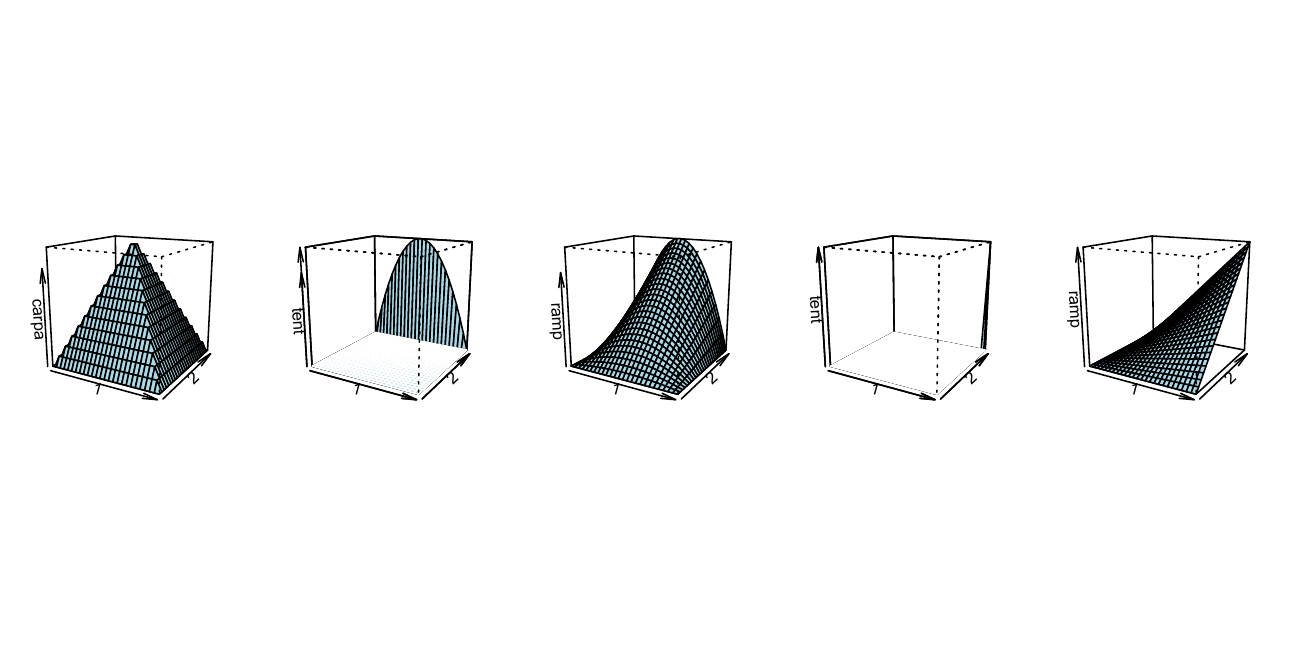}}
\vspace{-18mm}
\centerline{\scriptsize $\quad\quad C_{\{1,2\}}$-tent $\quad\quad\quad$ $C_{\{1\}}$-tent and its associated ramp $\quad\quad$ $C_{\{\emptyset\}}$-tent and its associated ramp}

\end{figure}

\medskip

The previous  notations allow us to state the following theorem:

\begin{theorem}\label{teorampa}$ $

\begin{itemize}
\item[(i)] If $g:C\rightarrow\R$ vanishes on $\partial^-C$, then there exists a unique decomposition
\begin{equation}\label{descg}g=\sum_{H\subset J}R_H\end{equation}
where each $R_H$ is an $H$-ramp associated to an $H$-tent $T_H$.
\item[(ii)] The maps $g\mapsto {\cal R}_Hg:=R_H$ and $g\mapsto {\cal T}_Hg:=T_H$ are linear, preserve the continuity and satisfy $\sup_{\sbm t\in C}{\cal R}_Hg(\bm t)\leq K_p\sup_{\sbm t\in C}g(\bm t)$, $\sup_{\sbm t\in C}{\cal T}_Hg(\bm t)\leq K_p\sup_{\sbm t\in C}g(\bm t)$, where $K_p$ depends on $p$ but not on $g$.
\end{itemize}\end{theorem}

The proof of this Theorem, where the ramps $R_H$ and tents $T_H$ are obtained constructively,  is deferred to the Appendix (\S\ref{appe}).

\subsection{Brownian tents, ramps, and a construction of Brownian sheet}\label{Brownianramps}
The $p$-parameter Wiener process or $p$-Brownian sheet on $C$ is the family of Gaussian centred variables $\{W(\bm t):\bm t\in C\}$ with covariances $\E W(\bm s)W(\bm t)=\ru\bm s\wedge\bm t\ru$ (see \cite{Chentsov1956,Walsh1986}). Since $W$ is a.s. 0 on the lower border $\partial^-C$ of $C$, Theorem \ref{teorampa} can be applied to conclude that there is a unique decomposition of $W$ as a sum of ramps:
\begin{equation}
\label{descwiener}W(\bm t)=\sum_{H\subset J}R_H(\bm t),\quad R_H(\bm t)=\ru\bm t_{J\setminus H}\ru T_H(\bm t_H).
\end{equation}
An important property of this decomposition is stated below in the Corollary of Theorem \ref{const}.

\begin{definition} $ $

\begin{itemize}
\item A Brownian $H$-tent is a centred Gaussian process on $C_H$ with covariances \begin{equation}\label{covtent}\E T_H(\bm s)T_H(\bm t)=\prod_{j\in H}(s_j\wedge t_j-s_jt_j).\end{equation}
\item A Brownian $H$-ramp is a $H$-ramp associated to a Brownian $H$-tent.\end{itemize} 
\end{definition}

Since  $\var T_H$  vanishes on  $\partial_H C_H$, the Brownian  $H$-tents are almost surely $H$-tents. These processes are already referred to as Brownian pillows in the mathematical literature (see, for instance, \cite{Koning&Protasov(2001),Hashova(2008)}). While an explicit definition of the Brownian pillow is provided in \cite{zhang2017}, the term ``pillow'' is also used in the same article to describe a multivariate extension associated with a different decomposition of the Brownian sheet. To avoid potential misunderstandings, we prefer to retain the name ``Brownian tent'' here.

\begin{theorem}\label{const} {\sc Construction of a $p$-Brownian sheet as a sum of independent Brownian ramps.} 

Let $\{T_H^*:H\subset J\}$ be a family of independent Brownian $H$-tents, and $R_H^*$ their corresponding Brownian ramps. 
Then the sum $W^*=\sum_{H\subset J}R_H^*$ is a $p$-Brownian sheet.\end{theorem}

{\bf Proof}. 
Because of the independence of the ramps, the covariances of the sum $W^*=\sum_{H\subset J}R_H^*$\ are
\begin{eqnarray*}\cov(W^*(\bm s),W^*(\bm t))&=&\sum_{H\subset J}\cov(R_H^*(\bm s)R_H^*(\bm t))\\&=&\sum_{H\subset J}\ru\bm s_{J\setminus H}\ru\ru\bm t_{J\setminus H}\ru\prod_{j \in H}(s_j\wedge t_j-e_jt_j)
\\&=&\ru\bm s\ru\ru\bm t\ru\sum_{H\subset J}\prod_{j\in H}\alpha_j \quad\mbox{ with }\alpha_j=\left(\frac{s_j\wedge t_j}{s_jt_j}-1\right).\end{eqnarray*}

On the other hand, a simple manipulation shows that the covariances $\ru\bm s\wedge \bm t\ru$ of the $p$-Brownian sheet 
are the same as the covariances of $W^*$: 
$$\ru\bm s\wedge \bm t\ru=\prod_{j=1}^ps_j\wedge t_j=\ru\bm s\ru\ru\bm t\ru\prod_{j=1}^p\frac{1}{s_j\vee t_j}=\ru\bm s\ru\ru\bm t\ru\prod_{j=1}^p(1+\alpha_j)=\ru\bm s\ru\ru\bm t\ru\sum_{H\subset J}\prod_{j\in H}\alpha_j$$ 
thus proving our statement.\qed

Since $W$ is a copy of $W^*$, the ramps ${\cal R}_HW$ are copies of the Brownian $H$-ramps $R_H^*$ and the tents ${\cal T}_HW$ are copies of the Brownian $H$-tents $T_H^*$, hence we have the following result.

\begin{corollary}\label{cor1}
The terms in the sum of equation (\ref{descwiener}) are independent Brownian ramps. 
\end{corollary}

\begin{remark} In these arguments the uniqueness of the decomposition and continuity of the maps ${\cal R}_H$ both established by Theorem \ref{teorampa} play a decisive role.\end{remark}

\begin{figure}
\caption{Decomposition of Brownian Sheet ($p=2$).}\label{figdec}
\centerline{\includegraphics[scale=0.60]{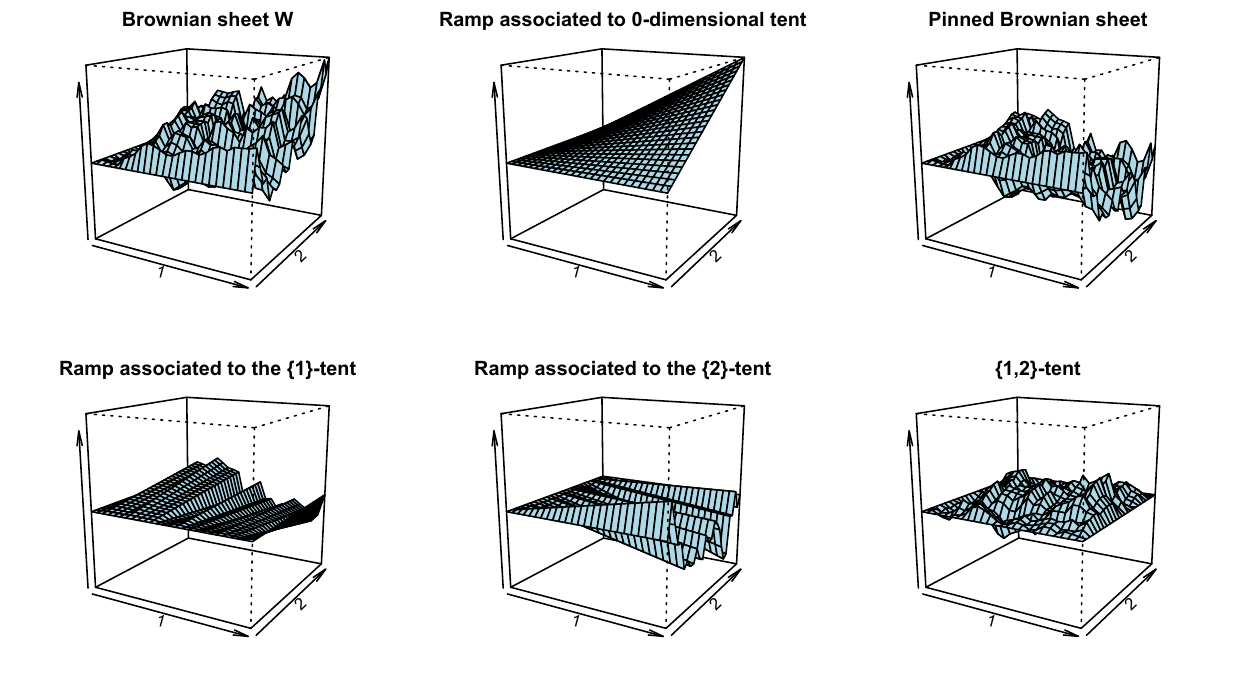}}

\end{figure}

The particular case of decomposition for $p=2$ is depicted in Figure \ref{figdec}.

\subsection{Probability laws of Brownian tents and their squared norms}

The probability distribution of the standard Brownian bridge is well known (see for instance \cite{durbin1973}) and may be found by obtaining the eigenfunctions $\{\psi_{\nu}:=\sqrt 2\sin(\nu\pi u\}_{\nu\in\sbm N:=\{1,2,3,\dots\}}$ of the covariance kernel $s\wedge t-st$, with eigenvalues $\lambda_{\nu}:=\frac{1}{\nu^2\pi^2}$. These eigenfunctions constitute an orthonormal basis of $L^2([0,1])$ with the Lebesgue measure. Then the law of $b$  is described by the Karhunen--Lo\`eve expansion $b(t)=\sum_{\nu\in\sbm N}\sqrt{\lambda_{\nu}}Z_{\nu}\Psi_{\nu}(t)$, where $\{Z_{\nu}:\nu\in\bm N\}$ are  i.i.d. standard Gaussian variables.

Since the covariance (\ref{covtent}) of a Brownian $H$-tent is a product of univariate kernels, the map $f(\bm t)\mapsto \int_{C_H}\prod_{j\in H}(s_j\wedge t_j-s_jt_j)f(\bm s)d\bm s_H$ has eigenfunctions $\psi_{\sbm \nu}(\bm t)=\prod_{j\in H}\lambda_{\nu_j}\psi_{\nu_j}(t_j)$, where $\bm \nu\in \N^{\#H}$ is a multi-index with components $\nu_j\in\N$ for each $j\in H$. These eigenfunctions are a complete orthonormal set on $L^2(C_H)$ with the Lebesgue measure, so that the Karhunen -- Lo\`eve expansion
\begin{equation}\label{KL}T_H(\bm t_H)=\sum_{\sbm\nu\in\N^{\#H}}\sqrt{\lambda_{\sbm\nu}}Z_{\sbm\nu}\psi_{\sbm\nu}(\bm t)=\sum_{\sbm\nu\in\N^{\#H}}\frac{Z_{\sbm\nu}}{\ru\bm\nu\ru\pi^{\#H}}\prod_{j\in H}\left(\sqrt 2\sin(\nu_j\pi t_j)\right)\end{equation} holds with $\{Z_{\sbm\nu}:\bm\nu\in\N^{\#H}\}$ i.i.d. standard Gaussian.

The probability law of the Brownian tent squared norm is obtained by integrating the square of (\ref{KL}) on $C_H$, and recalling that the eigenfunctions are an orthonormal set one gets \begin{equation}\label{nortent}\|T_H\|^2=\sum_{\sbm\nu\in\sN^{\#H}}\frac{Z_{\sbm\nu}^2}{\ru\bm\nu\ru^2\pi^{2\#H}}.\end{equation} Its c.d.f. shall be denoted $P_H(x)=\pr\{\|T_H\|^2\leq x\}$. 

\bigskip

\section{Uniformity tests based on the $H$-tents of the empirical process} \label{uniftests}

\subsection{Two consistent tests based on the asymptotic laws of the $H$-tents}\label{astests}

Let ${\cal U}_n\!=\!\{\bm U_1,\bm U_2,\dots,\bm U_n\}$ be
 a sample of i.i.d. $C$-valued random variables \linebreak$\bm U_i$ $=(U_{i,1}, U_{i,2}, \dots, U_{i,p})\tras$, $i=1,2,\dots,n$,
  with continuous distribution $F$, \linebreak$F_{{\cal U}_n}(\bm t)$ $=\frac{1}{n}\sum_{i=1}^n\un_{\{\sbm U_i\leq\sbm t\}}$
 its empirical distribution function and 
 $W_{{\cal U}_n}(\bm t)$ $=$ \linebreak$\frac{1}{\sqrt n}\sum_{i=1}^n(\bm 1_{\{\sbm U_i\leq\sbm t\}}-\ru\bm t\ru)$ 
the empirical process with respect to the uniform distribution. 

Then it is well known that
\begin{itemize} \item If $F$ is uniform in $C$ , $W_{{\cal U}_n}$ converges in law to the pinned Brownian sheet 
$W_0(\bm t)=W(\bm t)-\ru\bm t\ru W(\bm 1)$ as $n$ goes to infinity (\cite{PYKE197545}), 
and because of the continuity of the map ${\cal T}$ defined in Theorem \ref{teorampa}-(ii) the $H$-tents $T_{n,H}:={\cal T}_HW_{{\cal U}_n}$ of $W_{{\cal U}_n}$ converge jointly in law to the $H$-tents of $W_0$, which are the $H$-tents $T_H:={\cal T}_HW$ for $H\not=\emptyset$. This implies that $\|T_{n,H}\|^2\stackrel{\cal D}{\rightarrow}\|T_H\|^2$ jointly for $H\subset J,H\not=\emptyset$, and hence\begin{itemize}\item the vector $\bm p_n$ of components 
\begin{equation}\label{pnH}p_{n,H}=1-P_H(\|T_{n,H}\|^2)\end{equation} is asymptotically uniform on $[0,1]^{2^p-1}$;
 \item consequently the random variables $Q(1-p_{n,H})$, where $Q$ is the quantile function of the squared standard normal are asymptotically i.i.d.$\sim\chi^2$ with one degree of freedom, so that $\bm S:=\sum_{H\subset J, H\not=\emptyset}Q(1-p_{n,H})$ converges in law to a $\chi^2_f$ distribution with $f:=2^p-1$ degrees of freedom.
\end{itemize}
\item If $F$ is not uniform in $C$, $\lim_{n\rightarrow\infty}\|W_{{\cal U}_n}\|^2=\infty$ a.s.,
 so that the triangle inequality applied to (\ref{descwiener}) implies that at least for one non-empty $H$,  $\lim_{n\rightarrow\infty}\|T_{n,H}\|^2=3^{\#J\setminus H}\lim_{n\rightarrow\infty} \|R_{n,H}\|^2=\infty$ 
and hence at least one component of $\bm p_n$ tends to zero and $\bm S$ tends to infinity a.s.
\end{itemize}

The aforementioned dichotomy supports the rejection of the null hypothesis that $F$ follows a uniform distribution when ${\bm p}_{\mbox{\scriptsize min}} := \min\{p_{n,H} : H \subset J, H \neq \emptyset\}$ is smaller than a given constant, and likewise when $\bm{S}$ exceeds a specific constant. Both procedures result in consistent tests.
We refer to the test with the rejection region $\bm{p}_{\mbox{\scriptsize min}} < c$ as the minimum asymptotic test, abbreviated as m-as-test, and the test with the rejection region $S > k$ as the sum asymptotic test, or s-as-test. For large values of $n$, the asymptotic significance level of the m-as-test is $\alpha = 1 - (1 - c)^{2^p - 1}$. Meanwhile, selecting $k =Q_{2^p-1}(1-\alpha)$--the $1-\alpha$ quantile of the $\chi^2_{2p-1}$ distribution--yields an s-as-test with an asymptotic significance level of $\alpha$.

\subsection{Finite samples tests}

The practical implementation of the tests in section \ref{astests} requires the computation of the probabilities $P_H$ and the empirical tents $T_{n,H}$. 
These latter statistics are easily calculated as we show in \S\ref{ele}. However, the same is not true for the probability $P_H$ because the nice, compact formula (\ref{nortent}) does not allow for a simple calculation.

For this reason, we will replace the asymptotic tests with tests for each $n$, using the distribution $P_{n,H}$ of the statistics $\|T_{n,H}\|^2$ instead of its asymptotic distribution. Likewise, the exact calculation of probabilities needed to compute the new statistics will be replaced by estimates based on Monte Carlo simulations.

In summary, the decision procedures we propose are the following:

\subsubsection{The minimum test (abbreviated m-test) for samples of size $n$}\label{m-test}

\medskip

\begin{enumerate}

\item Generate the list $\cal H$ of nonempty subsets $H\subset J$,

\item for each $H\in{\cal H}$ compute the statistic $\|T_{n,H}\|^2$,

\item introduce the $p$-values $\tilde p_{n,H}=1-P_{n,H}(\|T_{n,H}\|^2)$, 

\item generate a large number $\{{\cal U}^r_n:r=1,2,\dots,R\}$ of independent samples of size $n$ of uniform random values on $C=[0,1]^p$,

\item  for each ${\cal U}^r_n$ compute the squared norms $\|T_H^r\|^2$ and estimate $\tilde p_{n,H}$ by means of the statistic $\hat p_{n,H}:=\frac{\sum_{r=1}^R\un_{\{\|T_H^r\|^2>\|T_H\|^2\}}+1}{R+1}$,

\item Reject the null hypothesis of uniformity if $\min_{\{H\in{\cal H}\}}\hat p_{n,H}$ is smaller than $1-(1-\alpha)^{1/\#{\cal H}}$ where $\alpha$ is the desired significance level of the test.

\end{enumerate}

\subsubsection{The sum test (abbreviated s-test) for samples of size $n$}\label{s-test}

\medskip

\begin{itemize}
\item[1.-5.] Repeat steps 1.-5. of \S\ref{m-test},

\item[6.] compute the statistic $\hat S=\sum_{H\in{\cal H}} Q_1(1-\hat p_{H})$ where $Q_1$ denotes the quantile function of the $\chi^2$ distribution with one degree of freedom,

\item[7.] reject the null hypothesis of uniformity if $\hat{\bm{S}}$ is greater than the quantile $1-\alpha$ of the $\chi^2$ distribution with $\#{\cal H}$ degrees of freedom.
\end{itemize}

\subsection{Partial tests}\label{partial}
If $p$ is very large, consistency can be sacrificed to reduce the number of empirical tents to be computed, by substituting a partial family of subsets of $J$, such as $\{H:H\subset J, 0<\#H\leq h\}$ with $h<p$ for the whole family $\cal H$ of nonempty subsets in the steps of sections \S\ref{m-test} and \S\ref{s-test}.

 Table \ref{lapartial} reports the empirical powers of those tests for 6-dimensional normal copula alternatives and critical regions with $h=1,2,3,4,5,6$.
 
 \medskip
 
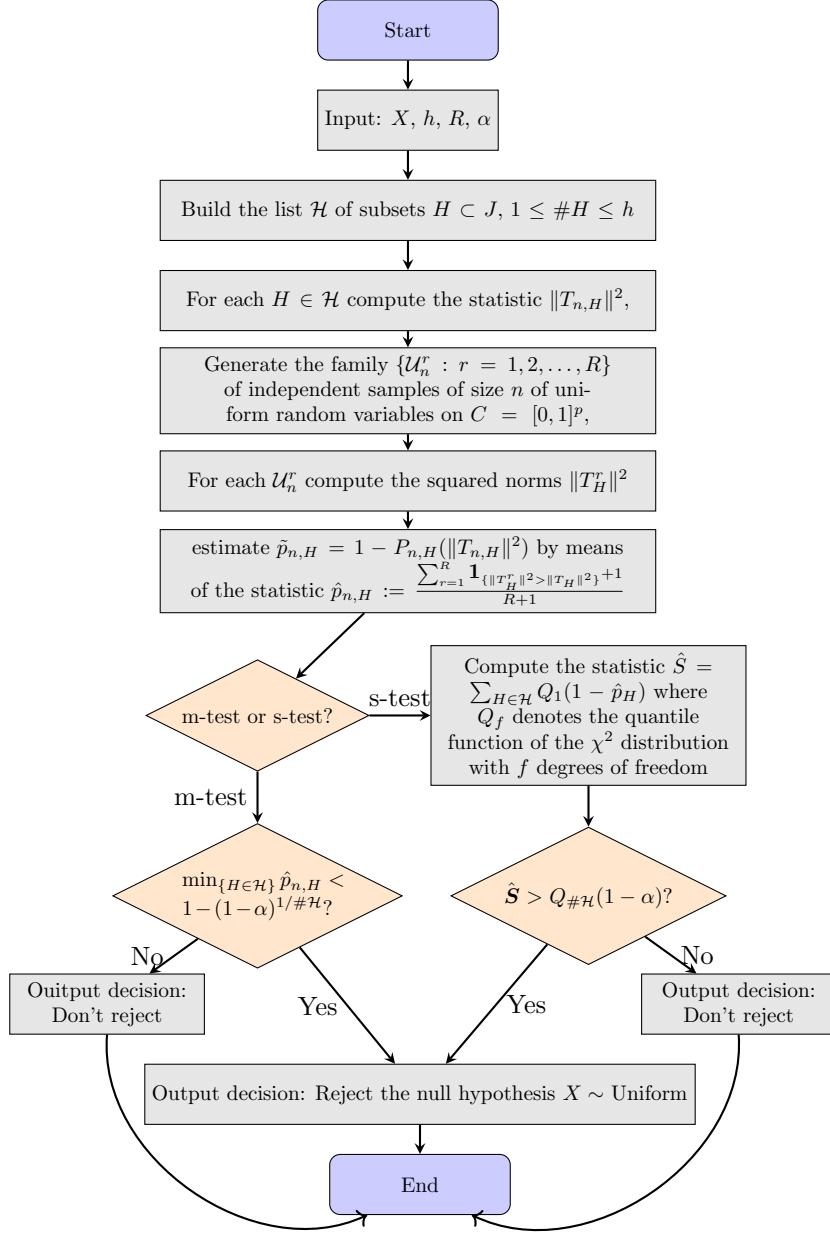
\begin{figure}\begin{center}

 \caption{A schematic summary of the proposed procedures to test uniformity of a multivariate sample $X$ with significance level $\alpha$} 
 \label{eldiag}
 
\begin{tikzpicture}[node distance=1.5cm]

\node[scale=.8] (start) [startstop] {Start};
\node[scale=.8] (input) [process, below of=start] {Input: $X$, $h$, $R$, $\alpha$};
\node[scale=.8] (buildH) [process, below of=input,text width=8cm] {Build the list $\cal H$ of subsets $H\subset J$, $1\leq \#H\leq h$};
\node[scale=.8] (compTH) [process, below of=buildH,text width=8cm] {For each $H\in{\cal H}$ compute the statistic $\|T_{n,H}\|^2$,};
\node[scale=.8] (GenU) [process, below of=compTH,text width=8cm] { Generate the family $\{{\cal U}^r_n:r=1,2,\dots,R\}$ of independent samples of size $n$ of uniform random variables on $C=[0,1]^p$,};
\node[scale=.8] (compTHr) [process, below of=GenU,text width=8cm] {For each ${\cal U}^r_n$ compute the squared norms $\|T_H^r\|^2$};
\node[scale=.8] (iestip) [process, below of=compTHr,text width=8cm] {estimate $\tilde p_{n,H}=1-P_{n,H}(\|T_{n,H}\|^2)$ by means of the statistic $\hat p_{n,H}:=\frac{\sum_{r=1}^R\un_{\{\|T_H^r\|^2>\|T_H\|^2\}}+1}{R+1}$};
\node[scale=.8] (quetest) [decision, below of=iestip, yshift=-0.9cm,xshift=-2.5cm] {m-test or s-test?};
\node[scale=.8] (compS) [process, right of=quetest, xshift=4cm, text width=5cm] {Compute the statistic $\hat S=\sum_{H\in{\cal H}} Q_1(1-\hat p_{H})$ where $Q_f$ denotes
 the quantile function of the $\chi^2$ distribution with $f$ degrees of freedom };
\node[scale=.8] (ispsmall) [decision, below of=quetest, yshift=-1.5cm,text width=2.5cm] {$\min_{\{H\in{\cal H}\}}\hat p_{n,H}<1-(1-\alpha)^{1/\#{\cal H}}$?};
\node[scale=.8](isSlarge)[decision,below of=compS, yshift=-1.5cm,text width=3cm]{$\hat{\bm{S}}>Q_{\#{\cal H}}(1-\alpha)$?};
\node[scale=.8] (output) [process, below of=ispsmall,xshift=2.7cm,yshift=-1.8cm] {Output decision: Reject the null hypothesis $X\sim$ Uniform};
\node[scale=.8] (end) [startstop, below of=output] {End};
\node[scale=.8](don'tm)[process,left of=ispsmall,text width=3cm,xshift=-1cm,yshift=-1.8cm]{Ouitput decision: Don't reject};
\node[scale=.8](don'ts)[process,right of=isSlarge,text width=3cm,xshift=1cm,yshift=-1.8cm]{Output decision: Don't reject};

\draw [arrow] (start) -- (input);
\draw [arrow] (input) -- (buildH);
\draw [arrow] (buildH) -- (compTH);
\draw [arrow] (compTH) -- (GenU);
\draw [arrow] (GenU) -- (compTHr);
\draw [arrow] (compTHr) -- (iestip);
\draw [arrow] (iestip) -- (quetest);
\draw [arrow] (quetest) -- node[above] {s-test} (compS);
\draw [arrow] (compS) -- (isSlarge);
\draw [arrow] (quetest) -- node[left] {m-test} (ispsmall);
\draw [arrow] (ispsmall) -- node[left]{No} (don'tm);
\draw [arrow] (ispsmall) -- node[left] {Yes}(output);
\draw[arrow] (isSlarge) --node[right] {No} (don'ts);
\draw [arrow] (isSlarge) -- node[right] {Yes}(output);
\draw [arrow] (output) -- (end);
\draw[thick,->, bend right=60] (don'tm) to (end);
\draw[thick,->, bend left=60] (don'ts) to (end);
\end{tikzpicture}

 \end{center}\end{figure}

 Figure \ref{eldiag} contains a schematic description of both m- and s-tests, including the partial versions.
 
\subsection{Computing the $H$-tent of $W_{{\cal U}_n}$}\label{ele}
The linearity of the map ${\cal T}$ allows to express the $H$-tents of samples of size $n$ in terms of $H$-tents of samples of size one.

The empirical distribution function of the sample of size one ${\cal U}_1=\{\bm U\}, \bm U=(U_1,U_2,\dots,U_p)\tras$, is $F_{{\cal U}_1}(\bm t)=\un_{\{\sbm U\leq\sbm t\}}=\prod_{j=1}^p\un_{\{U_j\leq t_j\}}$ and the empirical process is $W_{{\cal U}_1}(\bm t)$ $=\un_{\{\sbm U\leq\sbm t\}}-\ru\bm t\ru$ $=\prod_{j=1}^p\un_{\{U_j\leq t_j\}}$ $-\prod_{j=1}^pt_j$.

The ramps and tents of $W_{{\cal U}_1}$ can be obtained inductively by applying the operations described in the proof of Theorem \ref{teorampa}, but a simpler alternative way is to develop
\begin{eqnarray*}\un_{\{\sbm U\leq\sbm t\}}&=&\prod_{j=1}^p\un_{\{U_j\leq t_j\}}=\ru\bm t\ru\prod_{j=1}^p\left(1+\left(\frac{\un_{\{U_j\leq t_j\}}}{t_j}-1\right)\right)\\&=&\ru\bm t\ru\sum_{H\subset J}\prod_{j\in H}\left(\frac{\un_{\{U_j\leq t_j\}}}{t_j}-1\right)
=\sum_{H\subset J}\ru\bm t_{J\setminus H}\ru\prod_{j\in H}(\un_{\{U_j\leq t_j\}}-t_j)\\&=&\ru\bm t\ru+\sum_{H\subset J, H\not=\emptyset}\ru\bm t_{J\setminus H}\ru\prod_{j\in H}(\un_{\{U_j\leq t_j\}}-t_j)\end{eqnarray*} so that \begin{equation}\label{des}W_{{\cal U}_1}(\bm t)=\sum_{H\subset J, H\not=\emptyset}\ru\bm t_{J\setminus H}\ru\prod_{j\in H}(\un_{\{U_j\leq t_j\}}-t_j).\end{equation}

The products $\prod_{j\in H}(\un_{\{U_j\leq t_j\}}-t_j)$ are $H$-tents because they vanish if $\bm t_H\in\partial_HC_H$, that is, if al least one of the $t_j$, for $j\in H$, is $0$ or $1$.
Therefore the right-hand term of (\ref{des}) is the decomposition of $W_{{\cal U}_1}$ as a sum of $H$-ramps, thus proving that $${\cal R}_H(W_{{\cal U}_1})=\ru\bm t_{J\setminus H}\ru\prod_{j\in H}(\un_{\{U_j\leq t_j\}}-t_j)\quad\mbox{ and }\quad{\cal T}_H(W_{{\cal U}_1})=\prod_{j\in H}(\un_{\{U_j\leq t_j\}}-t_j).$$

The first part of next statement follows by applying the linear operator ${\cal T}_H$ to the equality $W_{{\cal U}_n}(\bm t)=\frac{1}{\sqrt n}\sum_{i=1}^nW_{\{\sbm U_i\}}(\bm t)$
and part (ii) follows by noticing that the integral in 
$\frac{1}{n}\int_{C_H}\sum_{h,i=1}^n\prod_{j\in H}(\un_{U_{h,j}\leq t_{j}}-t_j)(\un_{U_{i,j}\leq t_{j}}-t_j)d\bm t$ conmutes not only with the sum but also with the product because of the factorization of the integrand:

\begin{theorem}\label{empiricaltheo}

\begin{itemize} \item[(i)] The $H$-tent of the empirical process is $$T_{n,H}(\bm t)=\frac{1}{\sqrt n}\sum_{i=1}^n\prod_{j\in H}(\un_{\{U_{i,j}\leq t_j\}}-t_j),$$ and
\item[(ii)] its squared norm is \begin{equation}\label{empsqnorm}\|T_{n,H}\|^2=\frac{1}{n}\sum_{h,i=1}^n\prod_{j\in H}\left(\frac{U_{h,j}^2+U_{i,j}^2}{2}-U_{h,j}\vee U_{i,j}+\frac{1}{3}\right).\end{equation}
\end{itemize}\end{theorem}

\section{A brief empirical description of powers}\label{compar}

\subsection{Empirical comparison of the powers of several uniformity tests}
 Mengta Yang and Reza Modarres (\cite{Yang2015}) compare the powers of their uniformity tests $Q_1, Q_2, Q_3$ based on the distances between observations and
  $C_N$ based on the norm of the transformation of the sample obtained by applying 
  elementwise the inverse of the standard Normal c.d.f. to the sample points,
   with the powers of tests $M^2$ based on three types of depth of the sample points (\cite{Hegazy1975}),
     ${BCV}$ based on their distances-to-boundary (\cite{berrendero2006testing}) and ${MST}$ based on a minimum covering tree of the sample graph
      (\cite{PETRIE2013253}).
 
 For that purpose they construct two tables, both for samples in $[0,1]^2$, one for alternatives of dependence and the other for alternatives of shape.
 
Their first table contains the power of the different uniformity tests against the following copula alternatives:
\begin{itemize}
\item[]   AMH (Ali--Mikhail--Haq): $C_{\theta}(u, v)=\frac{uv}{ 1-\theta(1-u)(1-v)}$, $-1 \leq \theta < 1$,
\item[] FGM (Farlie--Gumbe--Morgenstern): $C_{\theta} (u, v) = uv + \theta uv(1-u)(1-v)$, $\theta\in [-1, 1]$,
\item[]  Clayton: $C_{\theta}(u, v) = \max[u^{-\theta} + v^{-\theta}-1,0]^{-1/\theta}$, $\theta\in [-1,\infty)$ y $\theta\not=0$,
 and
\item[]  Plackett: $C_{\theta}(u, v) = uv$ if $\theta=0$, and for $\theta >0$, $C_{\theta}(u,v)=\frac{1}{2}(\theta-1)A-\sqrt{A^2-4uv(\theta-1)}$, with $A = 1+(u+v)(\theta-1)$.
\end{itemize}

Observe that AMH and Clayton are Archimedean copulas, while FGM and Placket are non Archimed\-ean. 
\medskip

 In their second table the alternatives are vectors with i.i.d. components distributed Beta$(\alpha,\beta)$ for several values of the parameters.  

 Our tables \ref{YMdep} and \ref{YMsha} reproduce  those of Yang and Modarres and add at the end of each line, the estimated powers of our m- and s-tests based on 1000 replications of samples for each size and alternative.

 \begin{table}[t]
 
\tiny
\caption{Empirical power of uniformity tests against copula alternatives}\label{YMdep}
\setlength{\tabcolsep}{3pt}
\renewcommand{\arraystretch}{1.2}
\begin{tabular}{c*{10}{p{6mm}}*{2}{p{8mm}}}
alternative&$n$&$M^2_S$&$M^2_L$&$M^2_T$&$C_N$&${BCV}$&${MST}$&Q1&Q2&Q3&m-test&s-test\\\hline
AMH $\theta=0.9$&10&0.376&0.038&0.062&0.056&0.056&0.066&0.065&0.121&0.127&0.137&0.144\\
&25&0.328&0.118&0.054&0.068&0.062&0.112&0.066&0.170&0.164 & {\bf 0.359}&0.339  \\
&50&0.504&0.166&0.060&0.078&0.072&0.154&0.063&0.233&0.204 & {\bf 0.695}&0.648\\\hline
FGM $\theta=1$&10&0.672&0.046&0.096&0.060&0.044&0.044&0.055&0.090&0.094 & 0.093&0.086 \\
&25&0.590&0.076&0.060&0.072&0.040&0.052&0.052&0.104&0.115 & 0.238 &0.250\\
&50&0.390&0.072&0.050&0.062&0.040&0.094&0.049&0.126&0.127 & {\bf 0.459}&0.431 \\\hline
Clayton $\theta=2$&10&0.384&0.016&0.078&0.088&0.078&0.164&0.097&0.257&0.237 & 0.372&0.319\\
&25&0.638&0.472&0.076&0.074&0.136&0.592&0.101&0.427&0.370& {\bf 0.888}&0.849 \\
&50&0.984&0.850&0.060&0.090&0.194&0.894&0.098&0.640&0.566 &{\bf  0.998}&0.998 \\\hline
Plackett $\theta=5$&10&0.572&0.026&0.064&0.078&0.051&0.082&0.078&0.162&0.153 & 0.185&0.171\\
&25&0.414&0.170&0.046&0.072&0.078&0.152&0.076&0.234&0.210 & {\bf 0.536}&0.513 \\
&50&0.632&0.356&0.038&0.086&0.082&0.270&0.071&0.349&0.295& {\bf 0.860}&0.839\\\hline
\end{tabular}

\vspace{2mm}

\scriptsize
The numbers in boldface point out the cases in which our m-test outperforms the others.

\end{table}

\vspace{-10mm}

\begin{table}[t]
\tiny
\caption{Empirical powers of uniformity tests against bivariate i.i.d. Beta alternatives}
\label{YMsha}
\setlength{\tabcolsep}{5pt}
\renewcommand{\arraystretch}{1.15}
\begin{tabular}{*{2}{p{1.6mm}}*{9}{p{6mm}}*{2}{p{8mm}}}
$\alpha$&$\beta$&$M^2_S$&$M^2_L$&$M^2_T$&$C_N$&${{BCV}}$&${MST}$&Q1&Q2&Q3&m-test&s-test\\\hline
.5&.5&{\it 0.140}&{\it 0.356}&{\it 0.472}&{\it 0.268}&0.998&{\it 0.106}&0.997&0.999&0.999& 0.444&0.683\\
&1&{\it 0.330}&{\it 0.242}&{\it 0.182}&1.000&{\it 0.976}&{\it 0.254}&{\it 0.184}&{\it 0.415}&{\it 0.386} & 0.998&1.000 \\
&2&{\it 0.950}&{\it 0.698}&{\it 0.090}&1.000&1.000&{\it 0.998}&{\it 0.998}&{\it 0.951}&{\it 0.991} & 1.000&1.000\\
&3&{\it 0.996}&{\it 0.776}&{\it 0.086}&1.000&1.000&1.000&1.000&1.000&1.000& 1.000&1.000\\ \hline
1&1&0.056&0.054&{\it 0.044}&0.056&0.056&{\it 0.030}&0.048&{\it 0.042}&0.074 & 0.048&0.048 \\
&2&{\it 0.124}&{\it  0.254}&{\it 0.018}&1.000&{\it 0.066}&{\it 0.856}&{\it 0.971}&{\it 0.495}&{\it 0.965}& 1.000&1.000 \\
&3&{\it 0.374}&{\it 0.456}&{\it 0.070}&1.000&{\it 0.426}&1.000&1.000&{\it 0.221}&1.000& 1.000&1.000\\\hline
2&2&0.262&0.222&{\it 0.070}&{\it 0.030}&0.992&0.880&1.000&0.949&1.000& 0.108&0.207\\
&3&{\it 0.172}&{\it 0.314}&{\it 0.096}&{\it 0.806}&0.998&0.998&1.000&1.000&1.000& 0.977&0.994  \\\hline
3&3&{\it 0.166}&{\it 0.426}&{\it 0.150}&{\it 0.030}&1.000&1.000&1.000&{\it 0.544} &1.000& 0.720&0.935  \\\hline
\end{tabular}

\vspace{2mm}
\scriptsize
Numbers in italics indicate the cases with power smaller than the power of our s-test.

\end{table}

\FloatBarrier

Our tests show a good performance in detecting copula alternatives, outperforming in some cases all the competitors. As for the alternatives with i.i.d. components, they exhibit a 
power similar to that of some of the others, occasionally surpassing some of them.

\subsection{Performance of our partial tests against six-dimensional copulas}

Table \ref{lapartial} shows the empirical powers of the partial tests proposed in \S\ref{partial}, applied to samples of 50 normal copulas distributed as $\Phi_R(\Phi^{-1}(U_1),$ $\Phi^{-1}(U_2),\dots,\Phi^{-1}(U_p))$, where $U_1,U_2,\dots,U_p$ are i.i.d. Uniform on $[0,1]$, $\Phi$ is the standard Normal c.d.f. and $\Phi_R$ is the c.d.f. of the centred Normal vector in $\R^p$ with variance\\
\centerline{$R=\left(\begin{array}{ccccc}
1&\rho&\rho&\dots&\rho\\
\rho&1&\rho&\dots&\rho\\
\rho&\rho&1&\dots&\rho\\
\multicolumn{5}{c}{\dotfill}\\
\rho&\rho&\rho&\dots&1\end{array}\right)$}

The powers of Yang and Modarres (YM) tests computed by using the R package SHT are added for the sake of comparison.

\begin{table}[h]
\caption{Empirical powers of the YM tests $C_N$, 
Q1, Q2, Q3 
and our partial m- and s-tests.
}\label{lapartial}

\scriptsize{
    \begin{minipage}{\textwidth}
        \centering

\begin{tabular}{cccccccc}\hline
&\multicolumn{4}{c}{Yang \& Modarres tests}&\multicolumn{2}{c}{partial tests}&\\
alternative&C$_N$&Q1&Q2&Q3&m-test&s-test&$h$\\ \hline
&&&&&0.041&0.050 &1\\
normal&&&&&0.065&{\bf 0.098}&2 \\
copula&0.052& 0.043& 0.055& 0.077& 0.023&{\bf 0.147}&3 \\
$\rho=0.05$&&&&& 0.000 &{\bf 0.193}&4\\
&&&&&0.000&{\bf 0.220 }&5\\
&&&&& 0.000&{\bf 0.228 }&6\\ \hline
&&&&& 0.037&0.053&1\\
normal&&&&& 0.113& {\bf 0.236}&2\\
copula&0.051& 0.050& 0.107& 0.113&   0.042&{\bf0.238}&3\\
$\rho=0.10$&&&&& 0.000&{\bf  0.286}&4\\
&&&&& 0.000&{\bf  0.306}&5\\
&&&&& 0.000&{\bf  0.308}&6\\ \hline
&&&&&0.034&0.056&1\\
normal&&&&&0.193&{\bf 0.477}&2\\
copula&0.056& 0.060& 0.220& 0.191& 0.079&{\bf 0.417}&3\\
$\rho=0.15$&&&&&0.000&{\bf 0.420}&4\\
&&&&&0.000&{\bf 0.426}&5\\
&&&&&0.000&{\bf 0.429}&6\\ \hline
&&&&& 0.036&0.057&1\\
normal&&&&&0.327&{\bf 0.706}&2\\
copula&0.058 &0.072& 0.366& 0.314&0.150&{\bf 0.603}&3\\
$\rho=0.20$&&&&&0.000&{\bf 0.601}&4\\
&&&&&0.000&{\bf 0.590}&5\\
&&&&&0.000&{\bf 0.588}&6\\ \hline
&&&&& 0.034&0.056&1\\
normal&&&&&0.684&{\bf 0.965}&2\\
copula&0.060 &0.091& 0.720& 0.683&0.414&{\bf 0.906}&3\\
$\rho=0.30$&&&&&0.000&{\bf 0.881}&4\\
&&&&&0.000&{\bf 0.861} &5\\
&&&&&0.000&{\bf 0.857} &6\\ \hline
&&&&& 0.030&0.063&1\\
normal&&&&&0.911&{\bf 1.000}&2 \\
copula&0.065& 0.116 &0.929& 0.954&0.744&{\bf 0.988}&3\\
$\rho=0.40$&&&&&0.000&{\bf 0.974 }&4\\
&&&&&0.000&{\bf 0.968 }&5\\ 
&&&&&0.000&{\bf 0.968} &6\\\hline
\end{tabular}

\vspace{2mm}

\scriptsize
The numbers in boldface point out the cases in which our s-test outperforms the others.

\end{minipage}
}

\footnotesize \hspace{-1.8mm}

\bigskip

\end{table}

\subsection{Performance of our partial tests against six-dimensional copulas}

Table \ref{lapartial} shows the empirical powers of the partial tests proposed in \S\ref{partial}, applied to samples of 50 normal copulas distributed as $\Phi_R(\Phi^{-1}(U_1),$ $\Phi^{-1}(U_2),\dots,\Phi^{-1}(U_p))$, where $U_1,U_2,\dots,U_p$ are i.i.d. Uniform on $[0,1]$, $\Phi$ is the standard Normal c.d.f. and $\Phi_R$ is the c.d.f. of the centred Normal vector in $\R^p$ with variance\\
\centerline{$R=\left(\begin{array}{ccccc}
1&\rho&\rho&\dots&\rho\\
\rho&1&\rho&\dots&\rho\\
\rho&\rho&1&\dots&\rho\\
\multicolumn{5}{c}{\dotfill}\\
\rho&\rho&\rho&\dots&1\end{array}\right)$}

The powers of Yang and Modarres (YM) tests computed by using the R package SHT are added for the sake of comparison.

\section{Final comments} Both m- and s-tests are competitive with other tests proposed in the statistical literature,
and show a good performance in detecting copula alternatives.
It should be noted that the simplicity of the formula (\ref{empsqnorm}) for obtaining the $\|T_H\|^2$ allows for a simple calculation of the test statistics.

Neither of the two new tests is more powerful than the other. Which has better results depends on the alternatives.

\bigskip

Acknowledgement.

The first author was supported by grant  PID2021-123733NB-I00, Ministerio de Ciencia, Innovaci\'on y Universidades, Spain.

\addvspace{10pt}

\section{Appendix}

\subsection*{Proof of Theorem \ref{teorampa}}\label{appe}

In order to prove the Theorem 
we apply the following Lemma:

\begin{lemma} \label{lemarampa}Each $H$-ramp vanishes on all faces $C_{H'}$ with $\#H'=\#H$, $H'\not=H$.
\end{lemma}

{\bf Proof of Lemma \ref{lemarampa}} Assume that $R_H$ is the ramp associated to the $H$-tent $T_H$. There exist $j\in H\setminus H'$, so that if $\bm t\in C_{H'}$, then $\bm t_j=1$ and consequently the $j$-th component of $\bm t_H$ is one, that is, $\bm t_H$ belongs to $\partial C_H$, this implies that $T_H(\bm t_H)=0$, and hence the $H$-ramp $R_H(\bm t)=\ru\bm t_{J\setminus H}\ru T_H(\bm t_H)$ vanishes.\qed

\medskip

{\bf Proof of Theorem \ref{teorampa}}  The map $T_{\emptyset}:\mathbf{1}\mapsto g(\mathbf{1})$ is a $\emptyset$-tent, with $\emptyset$-ramp $R_{\emptyset}(\bm t)=\ru\bm t\ru g(\mathbf{1})$. 
We introduce now the function $g_0=g-R_{\emptyset}$ that vanishes on $C_{\emptyset}$.

For each $H$ with cardinal one, the restriction of $g_0$ to $C_H$ is a tent, because $\bm t\in \partial_H C_H$ implies $\bm t\in C_{\emptyset}$ or
 $\bm t\in\partial^-C$. Its ramp $R_H$ vanishes on the other faces of dimension $\#H$, as stated in Lemma \ref{lemarampa}, so that the function $g_1=g_0-\sum_{\#H=1}R_H=g-\sum_{\#H\leq 1}R_H$ vanishes on $\{C_H:\#H\leq 1\}$.

Now set $h=2,3,\dots,p$ succesively, and 
\begin{itemize}\item recall that $g_{h-1}$ vanishes on $\{C_H:\#H<h\}$, and notice that for $\#H=h$, the restriction of $g_{h-1}$ to $C_H$ is a tent, because $\bm t\in\partial_HC_H$ implies that either $\bm t$ belongs to a face of dimension $h-1$ or it belongs to $\partial^-C$,
 \item let $R_H$ denote the ramp associated to $g_{h-1}(t_H)$, that vanishes on the remaining faces of dimension $h$ as implied by Lemma \ref{lemarampa}, 
 \item introduce the function $g_h=g_{h-1}-\sum_{\#H=hj}R_H=g-\sum_{\#H\leq h}R_H$ that vanishes on $C_H$ for all $H$ with cardinal smaller or equal than $h$.
\end{itemize}

In particular, for $h=p$, $g_p=g-\sum_{\#H\leq p}R_H=0$ is the same as (\ref{descg}) thus proving the statement (i) of the Theorem.

The ramps $R_H$ are obtained by applying stepwise to $g$ two kind of operations:\begin{enumerate}\item restrictions of the domain and multiplication by continuous factors smaller or equal than one, which are linear, preserve the continuity and do not increase the supremum of the absolute value, and \item subtracting from $g$ $p\choose h$ ramps obtained by the first kind of operations, which preserves the continuity and the linearity.
\end{enumerate}
Let us denote $\rho_h=\sup_{\sbm t\in C}|g_{h-1}|$, so that for $\#H=h$, $\sup_{\sbm t\in C}|R_H|\leq \rho_h$ and therefore $\rho_{h+1}\leq (1+{p\choose h}) \rho_h$.
If $\rho=\sup_{\sbm t\in C}|g(\bm t)|$ then $\rho_1\leq (1+{p\choose 0})\rho=2\rho$, $\rho_2\leq(1+{p\choose 1})\rho_1=2(p+1)\rho$ and, in general $\rho_h\leq \prod_{j=0}^{h-1}(1+{p\choose j})\rho$. Therefore, the statement (ii) holds with $K_p=\prod_{j=0}^{p-1}\left(1+{p\choose j}\right)$.\qed

\bibliography{multinormality.bib}

\begin{thebibliography}{10}

\bibitem{PETRIE2013253}
An empirical study of tests for uniformity in multidimensional data.
\newblock {\em Computational Statistics \& Data Analysis}, 64:253--268, 2013.

\bibitem{berrendero2006testing}
José~R Berrendero, Antonio Cuevas, and Francisco V\'azquez-grande.
\newblock Testing multivariate uniformity: The distance-to-boundary method.
\newblock {\em Canadian Journal of Statistics}, 34(4):693--707, 2006.

\bibitem{Chentsov1956}
N.~N. Chentsov.
\newblock Weiner random fields depending on several parameters.
\newblock {\em Dokl. Akad. Nauk SSSR (N.S.)}, pages 607--609, 1956.

\bibitem{durbin1973}
J.~Durbin.
\newblock {\em {Distribution Theory for Tests Based on the Sample Distribution
  Function}}.
\newblock SIAM, 1973.

\bibitem{ebner2020testing}
Bruno Ebner, Franz Nestmann, and Matthias Schulte.
\newblock Testing multivariate uniformity based on random geometric graphs.

\bibitem{Hashova(2008)}
Enkelejd Hashorva.
\newblock Boundary non-crossings of brownian pillow.
\newblock {\em Journal of Theoretical Probability}, 23, 09 2008.

\bibitem{Hegazy1975}
Y.A.S. Hegazy and J.R. Green.
\newblock Some new goodness-of-fit tests using order statistics.
\newblock {\em JRStat Soc}, 24:299–308, 1975.

\bibitem{Koning&Protasov(2001)}
Alex Koning and Vladimir Protasov.
\newblock Tail behaviour of gaussian processes with applications to the
  brownian pillow.
\newblock {\em Journal of Multivariate Analysis}, 87:370--397, 01 2001.

\bibitem{PYKE197545}
Ronald Pyke.
\newblock Multidimensional empirical processes: Some comments.
\newblock In Madan~Lal Puri, editor, {\em Statistical Inference and Related
  Topics}, pages 45--58. Academic Press, 1975.

\bibitem{rosenblatt1952}
M.~Rosenblatt.
\newblock Remarks on a multivariate transformation.
\newblock {\em The Annals of Mathematical Statistics}, 23(3):470–472., 1952.
\newblock http://www.jstor.org/stable/2236692.

\bibitem{Walsh1986}
John~B WALSH.
\newblock {\em An introduction to stochastic partial differential equations}.
\newblock Springer, 1986.

\bibitem{Yang2015}
Mengta Yang and Reza Modarres.
\newblock Multivariate tests of uniformity.
\newblock {\em Stat Papers}, 58:627,639, 2015.

\bibitem{zhang2017}
Tonglin Zhang.
\newblock {On Independence and Separability between Points and Marks of Marked
  Point Processes, Supplementary material}.
\newblock {\em Statistica Sinica}, 27(1):109--114, 2017.

\end{thebibliography}

\end{document}